\documentclass[submission%
]{arxiv}

\usepackage[latin1]{inputenc}
\usepackage{subfigure}
\usepackage{amssymb,amsfonts,amsmath}
\usepackage{paralist}
\usepackage{pstricks}
\usepackage{url}

\usepackage[round]{natbib}

\author{David Einstein\thanks{Email: \email{deinst@gmail.com}.}
\and James Propp\thanks{Homepage: \url{http://jamespropp.org}.}}
\title[Piecewise-linear and birational toggling]
{Piecewise-linear and birational toggling}
\address{Department of Mathematical Sciences, UMass Lowell, USA}

\begin{document}
\maketitle
\begin{abstract}
\paragraph{Abstract.}
We define piecewise-linear and birational analogues of
the toggle-involutions on order ideals of posets
studied by Striker and Williams and use them 
to define corresponding analogues of rowmotion and promotion
that share many of the properties of combinatorial rowmotion and promotion.
Piecewise-linear rowmotion (like birational rowmotion)
admits an alternative definition
related to Stanley's transfer map for the order polytope;
piecewise-linear promotion relates to Sch{\"u}tzenberger promotion
for semistandard Young tableaux.
The three settings for these dynamical systems
(combinatorial, piecewise-linear, and birational)
are intimately related:
the piecewise-linear operations arise as tropicalizations
of the birational operations,
and the combinatorial operations arise as restrictions
of the piecewise-linear operations
to the vertex-set of the order polytope.
In the case where the poset is of the form $[a] \times [b]$,
we exploit a reciprocal symmetry property
recently proved by Grinberg and Roby
to show that birational rowmotion 
(and consequently piecewise-linear rowmotion) is of order $a+b$.
This yields a new proof of a theorem of Cameron and Fon-der-Flaass.
Our proofs make use of the correspondence 
between rowmotion and promotion orbits
discovered by Striker and Williams, which we make more concrete.
We also prove some homomesy results,
showing that for certain functions $f$,
the average value of $f$ over each rowmotion/promotion orbit
is independent of the orbit chosen.

\end{abstract}

\noindent
{\sc Note:} This is essentially a synopsis of the longer article-in-progress
\cite{einsteinpropp}.  It was prepared for FPSAC 2014, and will appear along
with the other FPSAC 2014 extended abstracts in a special issue of the journal
Discrete Mathematics and Theoretical Computer Science.

\newtheorem{theorem}{Theorem}
\newtheorem{lemma}[theorem]{Lemma}
\newtheorem{corollary}[theorem]{Corollary}
\newtheorem{proposition}[theorem]{Proposition}
\newtheorem{fact}[theorem]{Fact}
\newtheorem{observation}[theorem]{Observation}
\newtheorem{claim}[theorem]{Claim}

\newtheorem{definition}[theorem]{Definition}
\newtheorem{example}[theorem]{Example}
\newtheorem{conjecture}[theorem]{Conjecture}
\newtheorem{open}[theorem]{Open Problem}
\newtheorem{problem}[theorem]{Problem}
\newtheorem{question}[theorem]{Question}

\newtheorem{remark}[theorem]{Remark}
\newtheorem{note}[theorem]{Note}

\newcommand{\cA}{{\mathcal{A}}}
\newcommand{\cB}{{\mathcal{B}}}
\newcommand{\cC}{{\mathcal{C}}}
\newcommand{\cF}{{\mathcal{F}}}
\newcommand{\cH}{{\mathcal{H}}}
\newcommand{\cK}{{\mathcal{K}}}
\newcommand{\cO}{{\mathcal{O}}}
\newcommand{\cP}{{\mathcal{P}}}
\newcommand{\cS}{{\mathcal{S}}}
\newcommand{\row}{\rho}
\newcommand{\rowP}{\rho_{\cP}}
\newcommand{\rowH}{\rho_{\cH}}
\newcommand{\rowB}{\rho_{\cB}}
\newcommand{\pro}{\pi}
\newcommand{\proP}{\pi_{\cP}}
\newcommand{\proH}{\pi_{\cH}}
\newcommand{\proB}{\pi_{\cB}}
\newcommand{\proS}{\pi_{\cS}}
\newcommand{\phiP}{\phi_{\cP}}
\newcommand{\phiH}{\phi_{\cH}}
\newcommand{\phiB}{\phi_{\cB}}
\newcommand{\C}{\complexes}
\newcommand{\N}{\naturals}
\newcommand{\Q}{\rationals}
\newcommand{\R}{\reals}
\newcommand{\Z}{\integers}
\newcommand{\Sup}{S_{\rm up}}
\newcommand{\Sdown}{S_{\rm down}}
\newcommand{\psum}{\mathbin{\|}}
\newcommand{\covers}{\gtrdot}
\newcommand{\coveredby}{\lessdot}
\newcommand{\symmdiff}{\bigtriangleup}
\newcommand{\bigsig}{{\sum}}
\newcommand{\polytope}{K}
\newcommand{\mR}{\stackrel{\rho_\cB}{\mapsto}}
\newcommand{\mP}{\stackrel{\pi_\cB}{\mapsto}}

\def\urltilde{\kern -.15em\lower .7ex\hbox{\~{}}\kern .04em}
\def\eset{\emptyset}
\def\bc#1#2{\left(\kern -2pt{#1\atop #2} \kern -2pt\right)}
\def\zz{\mathbb}
\def\ds{\displaystyle}

%
%

\begin{section}{Background}
\label{sec:background}

We assume readers are familiar with 
the definition of a finite poset $(P, \leq)$,
as for instance given in Ch.\ 3 of~\cite{stanley2011enumerative}.
Much of our work 
involves the case $P = [a] \times [b] = 
\{(i,j) \in \N \times \N: \ 1 \leq i \leq a, \ 1 \leq j \leq b\}$
with ordering defined by $(i,j) \leq (i',j')$ 
iff $i \leq i'$ and $j \leq j'$.  We put $n=a+b$. 

We write $x \coveredby y$ (``$x$ is covered by $y$'')
or equivalently $y \covers x$ (``$y$ covers $x$'') 
when $x < y$ and no $z \in P$ satisfies $x < z < y$.
We say $P$ is {\em ranked} if there is a function
$r: P \rightarrow \{0,1,2,\dots\}$
so that all minimal elements of $P$ have rank 0
and such that $x \coveredby y$ implies $r(x) = r(y)-1$.

An {\em rc-embedding} of a poset $P$ is 
defined by~\cite{prorow} as
a map $\pi: P \rightarrow \Z \times \Z$
such that $x$ covers $y$ iff
$\pi(x)-\pi(y)$ is $(1,1)$ or $(-1,1)$.
This yields a Hasse diagram for $P$ in which
all covering relations are edges of slope $\pm 1$.
In the case $P = [a] \times [b]$, we will adopt
the rc-embedding $\pi$ that sends $(i,j) \in P$ to $(j-i,i+j-2) \in \Z^2$.
The {\em ranks} 
(or, in the terminology of~\cite{prorow}, {\em rows})
are the subsets of $P$ that consist of
all $x \in P$ at a given height, or vertical position,
relative to the rc-embedding.
We define the {\em files} 
(or, in the terminology of~\cite{prorow}, {\em columns})
as the subsets of $P$ that consist of
all $x \in P$ at a given horizontal position relative to the rc-embedding.
For example, let $P = [2] \times [2]$,
and write $(1,1),(2,1),(1,2),(2,2)$ as $w,x,y,z$ for short,
with $w < x < z$ and $w < y < z$.
Our rc-embedding has
$\pi(w) = (0,0)$, $\pi(x) = (-1,1)$, $\pi(y) = (1,1)$, and $\pi(z) = (0,2)$.
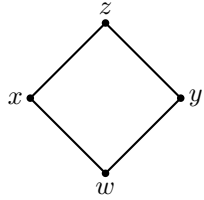
\begin{center}
\begin{pspicture}(-1,-.5)(1,2.5)
\psline(0,0)(-1,1)
\psline(0,0)(1,1)
\psline(-1,1)(0,2)
\psline(1,1)(0,2)
\pscircle[fillstyle=solid,fillcolor=black](0,0){.05}
\pscircle[fillstyle=solid,fillcolor=black](0,2){.05}
\pscircle[fillstyle=solid,fillcolor=black](-1,1){.05}
\pscircle[fillstyle=solid,fillcolor=black](1,1){.05}
\rput(0,-.2){$w$}
\rput(0,2.2){$z$}
\rput(-1.2,1){$x$}
\rput(1.2,1){$y$}
\end{pspicture}
\end{center}

The ranks are $\{w\}$, $\{x,y\}$, and $\{z\}$,
and the files are $\{x\}$, $\{w,z\}$, and $\{y\}$.
We number the ranks of $[a] \times [b]$ from 0 (bottom) to $n-1$ (top), 
and we number the files of $[a] \times [b]$ from 1 (left) to $n$ (right).
That is, for $P = [a] \times [b]$, $(i,j) \in P$ belongs to
the $(i+j-2)$nd rank of $P$ and 
the $(j-i+a)$th\footnote{Note that 
$j-i+a$ ranges from $1$ to $a+b-1=n-1$; 
this is slightly different from 
the indexing in~\cite{propproby}.} file of $P$.

We call $S \subseteq P$ an {\em order ideal} (or {\em downset}) of $P$
when $x \in S$ and $y \leq x$ imply $y \in S$.
We call $S \subseteq P$ a {\em filter} 
(or {\em upset}) of $P$
when $x \in S$ and $y \geq x$ imply $y \in S$.
We call $S \subseteq P$ an {\em antichain}
when $x, y \in S$ and $x \neq y$
imply that $x$ and $y$ are {\em incomparable}
(i.e., neither $x \leq y$ nor $y \leq x$). 
The sets consisting of the order ideals, filters, and antichains of $P$ 
are respectively denoted by $J(P)$, $\cF(P)$, and $\cA(P)$.

There are natural bijections 
$\alpha_1: J(P) \rightarrow \cF(P)$, 
$\alpha_2: \cF(P) \rightarrow \cA(P)$,
and $\alpha_3: \cA(P) \rightarrow J(P)$
given by the following recipes:
\begin{compactenum}
\item[(1)] for $I \in J(P)$, let $\alpha_1(I)$ be the complement $P \setminus I$;
\item[(2)] for $F \in \cF(P)$, let $\alpha_2(F)$ be the set of minimal elements of $F$
(i.e., the set of $x \in F$ such that $y < x$ implies $y \not\in F$); and
\item[(3)] for $A \in \cA(P)$, let $\alpha_3(A)$ be the downward saturation of $A$
(i.e., the set of $y \in P$ such that $y \leq x$ for some $x \in A$).
\end{compactenum}
The composition 
$\row := \alpha_3 \circ \alpha_2 \circ \alpha_1: J(P) \rightarrow J(P)$
is not the identity map;
e.g., it sends the full order ideal $I=P$
to the empty order ideal $I=\eset$.
(Note that~\cite{brouwerschrijver}
studied the closely related map
$F = \alpha_2 \circ \alpha_1 \circ \alpha_3: \cA(P) \rightarrow \cA(P)$.)

\cite{cameron1995orbits} gave an alternative characterization of $\row$.
Given $x \in P$ and $I \in J(P)$, let $\tau_x(I)$ 
(``$I$ toggled at $x$'' in Striker and Williams' terminology)
denote the set $I \symmdiff \{x\}$ if this set is in $J(P)$ and $I$ otherwise.
Equivalently, $\tau_x(I)$ is $I$
unless $y \in I$ for all $y \coveredby x$
and $y \not\in I$ for all $y \covers x$,
in which case $\tau_x(I)$ is $I \symmdiff \{x\}$.
(We will sometimes say that ``toggling $x$ turns $I$ into $\tau_x(I)$''.)
Clearly $\tau_x$ is an involution.
It is easy to show that $\tau_x$ and $\tau_y$ commute
unless $x \coveredby y$ or $x \covers y$.
If $x_1,x_2,\dots,x_{|P|}$ is any {\em linear extension} of $P$
(that is, a listing of the elements of $P$ such that
$x_i < x_j$ in $P$ implies $i<j$ in $\naturals$), then the composition 
$\tau_{x_1} \circ \tau_{x_2} \circ \cdots \circ \tau_{x_{|P|}}$
coincides with $\row$.
In the case where the poset $P$ is {\em ranked},
one natural way to linearly extend $P$ is to list the elements 
in order of increasing rank.
Given the right-to-left order of composition of
$\tau_{x_1} \circ \tau_{x_2} \circ \cdots \circ \tau_{x_{|P|}}$,
this corresponds to toggling the top rank first,
then the next-to-top rank, and so on, lastly toggling the bottom rank.
When $x$ and $y$ belong to the same rank of $P$,
the toggle operations $\tau_x$ and $\tau_y$ commute,
so even without using the theorem of Cameron and Fon-der-Flaass,
we can see that this composite operation on $J(P)$ is well-defined.
\cite{prorow} use the term ``row'' as a synonym for ``rank'',
and they refer to $\row$ as {\em rowmotion}.

In the example above,
under the action of $\alpha_1$, $\alpha_2$, and $\alpha_3$,
the order ideal $\{w,x\} \in J(P)$ gets successively mapped to
$\{y,z\} \in \cF(P)$, $\{y\} \in \cA(P)$, and $\{w,y\} \in J(P)$.
Under the action of $\tau_z$, $\tau_y$, $\tau_x$, and $\tau_w$,
the order ideal $\{w,x\} \in J(P)$ gets successively mapped to
$\{w,x\}$, $\{w,x,y\}$, $\{w,y\}$, and $\{w,y\}$ (all in $J(P)$).
In both cases we obtain $\row(\{w,x\}) = \{w,y\}$.

Note that if $x$ and $y$ belong to the same file,
the toggle operations $\tau_x$ and $\tau_y$ commute,
since neither of $x,y$ can cover the other.
Thus the composite operation of toggling the elements of $P$
from left to right is well-defined;
\cite{prorow} call this operation {\em promotion},
and show that it is conjugate to rowmotion in the toggle group
(the group generated by the toggle involutions). 
We denote this map by $\pro$.

\end{section}

\begin{section}{Piecewise-linear toggling}
\label{sec:PL}

Given a poset $P = \{x_1,\dots,x_p\}$ (with $p=|P|$)
and an rc-embedding of $P$,
let $\R^P$ denote the set of functions $f : P \rightarrow \R$;
we can represent such an $f$
as a {\em $P$-array} (or {\em array} for short)
in which the values of $f(x)$ for all $x \in P$ are arranged on the page
according to the rc-embedding of $P$ in the plane.
We will sometimes identify $\R^P$ with $\R^{p}$,
associating $f \in \R^P$ with $v = (f(x_1),\dots,f(x_p))$,
though this depends on the chosen ordering of the elements of $P$.
Let $\widehat{P}$ denote the augmented poset obtained from $P$ 
by adding two extra elements $\widehat{0}$ and $\widehat{1}$ 
(which we sometimes denote by $x_0$ and $x_{p+1}$)
satisfying $\widehat{0} < x < \widehat{1}$ for all $x \in P$. 
The {\em order polytope} $\cO(P) \subset \R^p$ 
(see~\cite{stanley86})
is the set of vectors $(\widehat{f}(x_1),\dots,\widehat{f}(x_p))$ in $\R^{p}$
arising from functions $\widehat{f} : \widehat{P} \rightarrow \R$ 
that satisfy $\widehat{f}(\widehat{0}) = 0$ and $\widehat{f}(\widehat{1}) = 1$ 
and are {\em order-preserving}
($x \leq y$ in $P$ implies $\widehat{f}(x) \leq \widehat{f}(y)$ in $\R$). 
In some cases it is better to work with the augmented vector 
$(\widehat{f}(x_0),\widehat{f}(x_1),\dots,\widehat{f}(x_p),\widehat{f}(x_{p+1}))$
in $\R^{p+2}$.
In either case we have a convex compact polytope.

For example, if $P = [2] \times [2] = \{w,x,y,z\}$,
then $\cO(P) = \{ v=(v_1,v_2,v_3,v_4) \in \R^4$ :
$0 \leq v_1$, $v_1 \leq v_2$, $v_1 \leq v_3$, 
$v_2 \leq v_4$, $v_3 \leq v_4$, and $v_4 \leq 1$\};
each such $v$ can be depicted as the $P$-array
\begin{displaymath}
\begin{array}{ccc}
    & v_4 &     \\[3mm]
v_2 &     & v_3 \\[3mm]
    & v_1 &   
\end{array}
\end{displaymath}
$\cO(P)$ is the convex hull of the vectors
$(0,0,0,0)$, $(0,0,0,1)$, $(0,0,1,1)$, 
$(0,1,0,1)$, $(0,1,1,1)$, and $(1,1,1,1)$,
which are precisely the vectors associated with the filters of $P$.
It is shown in~\cite{stanley86} that for any poset $P$,
the vertices of $\cO(P)$ correspond to the indicator functions
of the filters of $P$.

Given a convex compact polytope $\polytope$ in $\R^p$
(we are only concerned with the case $\polytope=\cO(P)$ here
but the definition makes sense more generally),
we define the {\em piecewise-linear toggle operation} 
$\tau_i$ ($1 \leq i \leq p$)
as the unique map from $\polytope$ to itself
whose action on the 1-dimensional cross-sections of $\polytope$ 
in the $i$th coordinate direction
is the linear map that switches the two endpoints of the cross-section.
That is, given $v = (v_1,\dots,v_p) \in \polytope$,
we define 
\begin{equation}
\label{eq:phiv}
\tau_i(v) = (v_1,\dots,v_{i-1},L+R-v_i,v_{i+1},\dots,v_p),
\end{equation}
where the real numbers $L$ and $R$
are respectively the left and right endpoints of the set 
$\{t \in \R: \ (v_1,\dots,v_{i-1},t,v_{i+1},\dots,v_p)
\in \polytope\}$,
which is a bounded interval because $K$ is 
convex and compact.\footnote{Note that $L$ and $R$ depend on
$v_1,\dots,v_{i-1},v_{i+1},\dots,v_p$,
though our notation suppresses this dependence.}
Since $L+R-(L+R-v_i)=v_i$, each toggle operation is an involution.

Similar involutions were studied by~\cite{kirillov1996groups}
in the context of Gelfand-Tsetlin triangles.
Indeed, one can view their action in our piecewise-linear toggling framework,
where instead of looking at the rectangle posets $[a] \times [b]$
one looks at the triangle posets
with elements $\{(i,j): \ 1 \leq i \leq j \leq N\}$
and covering-relations $(i,j-1) \coveredby (i,j)$ (for $1 \leq i \leq j \leq N$)
and $(i+1,j+1) \coveredby (i,j)$ (for $1 \leq i \leq j \leq N-1$).
Their ``elementary transformations'' (Definition 0.1) are our ``toggles''.

In the case where $\polytope$ is the order polytope of $P$
and a particular element $x \in P$ has been indexed as $x_i$,
we write $\tau_i$ as $\tau_x$. 
The $L$ and $R$ that appear in (\ref{eq:phiv})
are given by
\begin{equation}
\label{eq:L}
L = \max \{v_j: \ 0 \leq j \leq p+1, \ x_j \coveredby x_i\}
\end{equation}
and
\begin{equation}
\label{eq:R}
R = \min \{v_j: \ 0 \leq j \leq p+1, \ x_j \covers x_i\} .
\end{equation}
(One also has 
$L = \max \{v_j : x_j < x_i\}$ 
and 
$R = \min \{v_j : x_j > x_i\}$, 
but the formulas (\ref{eq:phiv})--(\ref{eq:R}) 
turn out to be the right ones to use 
when extending the operations $\tau_i$
from $\cO(P)$ to all of $\R^p$,
as well as the right ones to use
when lifting toggling to the birational setting
as described in the next section.)
It is easy to show that $\tau_x$ and $\tau_y$ commute
unless $x \coveredby y$ or $x \covers y$.
These piecewise-linear toggle operations $\tau_x$ are analogous to
the combinatorial toggle operations $\tau_x$ 
(and indeed the former generalize the latter
in a sense to be made precise below),
so it is natural to define 
piecewise-linear rowmotion $\rowP: \cO(P) \rightarrow \cO(P)$
as the composite operation accomplished by toggling from top to bottom
(much as ordinary rowmotion $\row: J(P) \rightarrow J(P)$
can be defined as the composite operation obtained by
toggling from top to bottom).  Likewise we can define
piecewise-linear promotion $\proP: \cO(P) \rightarrow \cO(P)$
as the composite operation accomplished by toggling from left to right.

Continuing the example $P = [2] \times [2] = \{w,x,y,z\}$
from section~\ref{sec:background},
let $v = (.1,.2,.3,.4) \in \cO(P)$.
Under the action of $\tau_z$, $\tau_y$, $\tau_x$, and $\tau_w$,
the vector $v$ gets successively mapped to
$(.1,.2,.3,.9)$, $(.1,.2,.7,.9)$, $(.1,.8,.7,.9)$, 
and $(.6,.8,.7,.9) = \rowP(v)$,
while under the action of $\tau_x$, $\tau_w$, $\tau_z$, and $\tau_y$,
the vector $v$ gets successively mapped to
$(.1,.3,.3,.4)$, $(.2,.3,.3,.4)$, $(.2,.3,.3,.9)$, 
and $(.2,.3,.8,.9) = \proP(v)$.

If $f$ is the indicator function of the filter $P \setminus I$,
then $\rowP(v)$ (resp.\ $\proP(v)$) is the indicator function
of the filter $P \setminus \row(I)$ (resp.\ $P \setminus \pro(I)$);
in this way $\rowP$ and $\proP$ generalize $\row$ and $\pro$.

In the full version of the article (\cite{einsteinpropp}),
we extend $\rowP$ and $\proP$ to all of $\R^p$, not just $\cO(P)$.
We also study a variant of these extended operations in which one takes 
$(\widehat{f}(\widehat{0}), \widehat{f}(\widehat{1})) = (0,0)$
instead of $(0,1)$;
although there is no longer an order polytope in the picture,
these ``homogeneous'' actions are easier to understand,
and capture most of the behavior of the general inhomogeneous case.

One can show that the action of the Sch{\"u}tzenberger promotion operator 
(which we denote by $\proS$) on the set of semistandard Young tableaux
of rectangular shape with $A$ rows and $B$ columns
having entries between 1 and $n$ is naturally conjugate to the action of 
the piecewise-linear promotion operator $\proP$
on the rational points in the order polytope of $P = [A] \times [n-A]$
with denominator dividing $B$.
(We are grateful to Alex Postnikov and Darij Grinberg
for explaining this to us.
For the original definition of promotion,
see~\cite{schutzenberger1972promotion};
for more modern treatments,
see~\cite{stanley2009promotion} and~\cite{vanleeuwen}.)
For example, take $A=2$, $B=3$, and $n=5$,
and consider the semistandard Young tableau
\begin{displaymath}
\begin{array}{ccc}
1 & 2 & 2 \\
3 & 5 & 5
\end{array}
\end{displaymath}
We represent the tableau $T$ as a Gelfand-Tsetlin triangle 
whose $i$th row ($1 \leq i \leq n$) lists, in decreasing order 
(with 0's appended or deleted from the end as needed),
the number of parts less than or equal to $n-i+1$
in the successive rows of the tableau:
\begin{displaymath}
\begin{array}{ccccccccc}
\overline{3} &   & \overline{3} &   & \underline{0} &   & \underline{0} &   & \underline{0} \\
  & \overline{3} &   & 1 &   & \underline{0} &   & \underline{0} &   \\
  &   & 3 &   & 1 &   & \underline{0} &   &   \\
  &   &   & 3 &   & 0 &   &   &   \\
  &   &   &   & 1 &   &   &   &
\end{array}
\end{displaymath}
This tableau splits into three parts: 
a triangle of $B$'s (overlined, with top row of length $A$),
a triangle of 0's (underlined, with top row of length $n-A$),
and an $(n-A)$-by-$(A)$ rectangle.
If we flip this rectangle across the line $x+y=0$,
so that the top corner becomes the left corner and vice versa,
we get a $P$-array with entries between 0 and $B$:
\begin{displaymath}
\begin{array}{cccc}
  & 3 &   &   \\
1 &   & 3 &   \\
  & 1 &   & 1 \\
  &   & 0 &
\end{array}
\end{displaymath}
If we divide each entry by $B$, we get a point $v(T)$ in $\cO(P)$
from which one can recover $T$ by reversing all the above steps.
One can show that $v(\proS(T))=\proP(v(T))$.
Indeed, the file-toggle operations
(in which one performs piecewise-linear toggling
at all $x \in P$ belonging to the $i$th file
of $[n-A] \times [A]$, with $1 \leq i \leq n$;
see~\ref{sec:filetoggle})
can be shown to correspond respectively to
the $n$ Bender-Knuth involutions on the Young tableau,
whose composition gives $\proS$.

The vertices of $\cO(P)$ correspond to
the 0,1-valued functions $f$ on $P$
with the property that $x \leq y$ in $P$ implies $f(x) \leq f(y)$ in $\{0,1\}$;
these are precisely the indicator functions of filters.
Filters are in bijection with order ideals
by way of the complementation map,
so the vertices of $\cO(P)$ are in bijection 
with the elements of the lattice $J(P)$.
Each toggle operation acts as a permutation on the vertices of $\cO(P)$. 
Indeed, if we think of each vertex of $\cO(P)$
as determining a cut of the poset $P$
into an upset (filter) $\Sup$ and a complementary downset (order ideal) $\Sdown$
(the pre-image of 1 and 0, respectively,
under the order-preserving map from $P$ to $\{0,1\}$),
then the effect of the toggle operation $\tau_x$ ($x \in P$)
is just to move $x$ from $\Sup$ to $\Sdown$ (if $x$ is in $\Sup$)
or from $\Sdown$ to $\Sup$ (if $x$ is in $\Sdown$)
unless this would violate the property
that $\Sup$ must remain an upset and $\Sdown$ must remain a downset.
In particular, we can see that when our point $v \in \cO(P)$ is a vertex
associated with the cut $(\Sup,\Sdown)$,
the effect of $\tau_x$ on $\Sdown$ is just
toggling the order ideal $\Sdown$ at the element $x \in P$.


\cite{cameron1995orbits} showed that
rowmotion acting on $J([a] \times [b])$ is of order $a+b$.
(Subsequently \cite{prorow} gave a simpler proof,
by showing that promotion is of order $a+b$
and that rowmotion is conjugate to promotion.)
The same is true of piecewise-linear rowmotion and promotion
acting on $\cO([a] \times [b])$: 

\begin{theorem}
\label{thm:pl-order}
For $P = [a] \times [b]$,
the maps $\rowP$ and $\proP$
are of order $a+b$.
\end{theorem}

It seems plausible that one might be able to deduce 
the order of $\rowP$ and $\proP$ from the order of $\row$ and $\pro$,
but we have not been able to find such an argument.\footnote{The 
Coxeter hyperplane arrangement of type $A$
divides the order polytope into simplices,
and on each simplex the maps $\rowP$ and $\proP$
are not just piecewise-linear but actually linear
(by which we really mean ``affine''),
and one might hope to base a proof of Theorem~\ref{thm:pl-order} on this;
unfortunately, the images of these simplices
under $\rowP$ and $\proP$ are not themselves
simplices in this dissection,
so the most simple sort of proof one might imagine
does not work.}
Instead, our proof of Theorem~\ref{thm:pl-order} 
detours through the notions of birational promotion and rowmotion.

\end{section}

\begin{section}{Birational toggling}
\label{sec:birational}

The definition of the piecewise-linear toggling operation
via formulas (\ref{eq:phiv})--(\ref{eq:R}) 
involves only addition, subtraction, min, and max.
Consequently one can define birational transformations on $(\R^+)^P$
with formal resemblance to the toggle operations on $\cO(P)$.
This transfer makes use of a dictionary in which 
0, addition, subtraction, max, and min
are respectively replaced by 
1, multiplication, division, addition, and parallel addition (defined below),
resulting in a subtraction-free rational expression.\footnote{The 
authors are indebted to Arkady Berenstein for pointing out
the details of this transfer of structure
from the piecewise-linear setting to the birational setting.}
Parallel addition can be expressed in terms of the other operations,
but taking a symmetrical view of the two forms of addition 
turns out to be fruitful.
Indeed, in setting up the correspondence we have a choice to make:
by ``series-parallel duality'',
one could equally well use a dictionary 
that switches the roles of addition and parallel addition.
We hope the choice that we have made here
will prove to be convenient.

For $x,y$ satisfying $x+y \neq 0$, we define 
the parallel sum of $x$ and $y$ as $x \psum y = xy/(x+y)$.
In the case where $x$, $y$ and $x+y$ are all nonzero,
$xy/(x+y)$ is equal to $1/(\frac1x+\frac1y)$,
which clarifies the choice of notation and terminology:
if two electrical resistors of resistance $x$ and $y$ are connected in parallel,
the compound circuit has an effective resistance of $x \psum y$.
If $x$ and $y$ are in $\R^+$, 
then $x+y$ and $x \psum y$ are in $\R^+$ as well.
Also, $\psum$ is commutative and associative,
so that a compound parallel sum $x \psum y \psum z \psum \cdots$
is well-defined; it equals the product $x y z \cdots$
divided by the sum of all products that omit exactly one of the variables,
and in the case where $x,y,z,\dots$ are all positive,
it can also be written as $1/(\frac1x+\frac1y+\frac1z+\cdots)$.

Given a non-empty set $S = \{s_1,s_2,\dots\}$,
let $\bigsig^+ S$ denote $s_1 + s_2 + \cdots$ 
and $\bigsig^{\psum} S$ denote $s_1 \psum s_2 \psum \cdots$.
Then for $v = (v_0,v_1,\dots,v_p,v_{p+1}) \in (\R^+)^{p+2}$ 
with $v_0 = v_{p+1} = 1$
and for $1 \leq i \leq p$ we define
\begin{equation}
\label{eq:phib}
\tau_i(v) = (v_0,v_1,\dots,v_{i-1},LR/v_i,v_{i+1},\dots,v_p,v_{p+1}),
\end{equation}
with
\begin{equation}
\label{eq:Lb}
L = \bigsig^{+} \{v_j: \ 0 \leq j \leq p+1, \ x_j \coveredby x_i\}
\end{equation}
and
\begin{equation}
\label{eq:Rb}
R = \bigsig^{\psum} \{v_j: \ 0 \leq j \leq p+1, \ x_j \covers x_i\}.
\end{equation}
We call the maps $\tau_i: (\R^+)^P \rightarrow (\R^+)^P$ 
given by (\ref{eq:phib})--(\ref{eq:Rb})
{\em birational toggle operations},
as opposed to the piecewise-linear toggle operations
treated in the previous section.\footnote{In principle
we should use a different symbol than $\tau_i$,
but in practice it should always be clear
whether we are referring to piecewise-linear operations
or birational operations.}
As the $0$th and $p+1$st coordinates of $v$ 
are not affected by any of the toggle operations,
we can just omit those coordinates,
reducing our toggle operations to actions on $(\R^+)^p$.
Since $LR/(LR/v_i) = v_i$,
each birational toggle operation is an involution
on the orthant $(\R^+)^{p}$.
As in the preceding section,
we identify $(\R^+)^p$ with $(\R^+)^P$.
The birational toggle operations are analogous to
the piecewise-linear toggle operations
(in a sense to be made precise below),
so it is natural to define 
{\em birational rowmotion} $\rowB: (\R^+)^P \rightarrow (\R^+)^P$
as the composite operation accomplished by toggling from top to bottom,
and to define
{\em birational promotion} $\proB: (\R^+)^P \rightarrow (\R^+)^P$
as the composite operation accomplished by toggling from left to right.

Continuing our running example $P = [2] \times [2] = \{w,x,y,z\}$,
let $v = (1,2,3,4) \in \R^P$,
corresponding to the positive function $f$
that maps $w,x,y,z$ to $1,2,3,4$, respectively,
with $f(\widehat{0}) = f(\widehat{1}) = 1$.
Under the action of $\tau_z$, $\tau_y$, $\tau_x$, and $\tau_w$,
the vector $v = (1,2,3,4)$ gets successively mapped to
$(1,2,3,\frac{5}{4})$,
$(1,2,\frac{5}{12},\frac{5}{4})$,
$(1,\frac{5}{8},\frac{5}{12},\frac{5}{4})$,
and $(\frac{1}{4},\frac{5}{8},\frac{5}{12},\frac{5}{4}) = \rowB(v)$.

For simplicity, we have defined
$\proB$ as a map from $(\R^+)^P$ to itself. 
However, $\proB$ can be extended to a map
from a dense open subset of $\R^P$ to itself,
and indeed, from a dense open subset $U$ of $\C^P$ to itself.
All expressions we consider are well-defined on the open orthant $(\R^+)^P$,
and all the theorems we prove amount to identities
that are valid when all variables lie in this orthant;
this implies that the identities hold
outside of some singular variety in $\C^P$.
Identifying the singular subvariety on which $\proB$ 
(or one of its powers) is undefined seems like an interesting question, 
but it is one that we leave to others.
Alternatively, Tom Roby has pointed out that one can replace $\R^+$ by
a ring of rational functions in formal indeterminates
indexed by the elements of $P$,
thereby avoiding the singularity issue
(once one checks that the rational functions in question
can be expressed as ratios of polynomials with positive coefficients).

Piecewise-linear rowmotion and promotion can be viewed
as tropicalizations of birational rowmotion and promotion.
To the extent that facts about birational toggling
can be formulated as (complicated but finite) identities 
in subtraction-free arithmetic,
the dictionary alluded to at the start of section~\ref{sec:birational}
allows one to carry the identities 
to the ``max, min, plus'' setting.\footnote{We
are indebted to Colin McQuillan and Will Sawin 
for clarifying this point; see~\cite{mathoverflow1}.}
For instance, when in a later section
we prove that $\rowB^n$ and $\proB^n$ act trivially on $(\R^+)^P$
(with $P=[a] \times [b]$ and $n=a+b$),
it will follow immediately
that $\rowP^n$ and $\proP^n$ act trivially on $\R^P$.
(Here we gloss over the role that $\widehat{0}$ and $\widehat{1}$ play.
Our treatment of birational toggling
assumes $\widehat{f}(\widehat{0}) = \widehat{f}(\widehat{1}) = 1$
but our treatment of piecewise-linear toggling
assumes $\widehat{f}(\widehat{0}) = 0 \neq 1 = \widehat{f}(\widehat{1})$.
The full version of the paper addresses this issue
with an appropriate dehomogenization lemma.)

\end{section}

\begin{section}{Birational rowmotion and Stanley's transfer map}

\label{sec:transfer}

Although most of our work with rowmotion
treats it as a composition of $|P|$ toggles
(from the top to the bottom of $P$),
we noted in section~\ref{sec:background}
that $\row$ can also be defined
as a composition of three operations
$\alpha_1$, $\alpha_2$, $\alpha_3$.\footnote{Indeed 
this was the way in which Brouwer and Schrijver
originally defined their operation $F$,
in the context of the Boolean lattices 
$[2] \times [2] \times \cdots \times [2]$.}
This alternative definition can be lifted 
to the piecewise-linear and birational settings.

For the piecewise-linear setting, we first recall 
the definition of the {\em chain polytope} $\cC(P)$ of a poset $P$
as defined by~\cite{stanley86}.
A {\em chain} in a poset $P$ is a totally ordered subset of $P$,
and a {\em maximal chain} in a poset $P$
is a chain that is not a proper subset of any other chain.
If the poset $P$ is ranked, with all maximal elements having the same rank,
then the maximal chains in $P$ 
are precisely those chains that contain an element of every rank.
The {\em chain polytope} of a poset $P$
is the set of maps from $P$ to $[0,1]$
such that for every chain $C$ in $P$
(or, equivalently, for every maximal chain $C$ in $P$),
\begin{equation}
\label{eq:chain}
\sum_{x \in C} f(x) \leq 1.
\end{equation}
Just as the vertices of the order polytope of $P$
correspond to the indicator functions of the filters of $P$,
the vertices of the chain polytope of $P$ 
correspond\footnote{One direction of this claim is easy:
since every antichain intersects every chain of $P$
in at most one element of $P$,
the indicator function of an antichain must correspond to a point in $\cC(P)$.
For the other direction,
see Theorem 2.2 of~\cite{stanley86}.}
to the indicator functions 
of the antichains of $P$.

Stanley defines the {\em transfer map} $\Phi: \cO(P) \rightarrow \R^P$
via the formula
\begin{equation}
\label{eq:transfer}
(\Phi f)(x) = \min \{ f(x)-f(y): y \in \widehat{P}, \ x \covers y \}
\end{equation}
for all $x \in P$ (recall that we have $f(\widehat{0}) = 0$).
Stanley proves that $\Phi$ is a bijection between $\cO(P)$ and $\cC(P)$
that carries the vertices of the former to the vertices of the latter.
The inverse of $\Phi$ is given by\footnote{This is not 
precisely the definition of $\Psi$ that Stanley gives,
but the two definitions are easily seen to be equivalent.}
\begin{equation}
\label{eq:inverse}
(\Psi g)(x) = \max \{ g(y_1)+g(y_2)+\cdots+g(y_k): \widehat{0} \coveredby y_1
\coveredby y_2 \coveredby \cdots \coveredby y_k = x \} .
\end{equation}

Let $\tilde{\cO}(P)$ be the set of order-reversing maps from $P$ to $[0,1]$.
We now define bijections 
$\alpha_1: \tilde{\cO}(P) \rightarrow \cO(P)$, 
$\alpha_2: \cO(P) \rightarrow \cC(P)$,
and $\alpha_3: \cC(P) \rightarrow \tilde{\cO}(P)$
given by the following recipes:
\begin{compactenum}
\item[(1)] for $f \in \tilde{\cO}(P)$, let $\alpha_1(f)$ be defined by 
\begin{displaymath} (\alpha_1(f))(x) = 1-f(x); \end{displaymath}
\item[(2)] for $f \in \cO(P)$, let $\alpha_2(f)$ be defined by
\begin{displaymath} (\alpha_2 f)(x) = 
\min \{ f(x)-f(y): y \in \widehat{P}, \ x \covers y \}; \end{displaymath}
and
\item[(3)] for $f \in \cC(P)$, let $\alpha_3(f)$ be defined by
\begin{equation}
\label{eq:longsum}
(\alpha_3 f)(x) = \max \{ f(y_1)+f(y_2)+\cdots+f(y_k): \ x = y_1
\coveredby y_2 \coveredby \cdots \coveredby y_k \coveredby \widehat{1} \} .
\end{equation}
\end{compactenum}
Note that $\alpha_2$ is $\Phi$
and that $\alpha_3$ is $\Psi$ (aka $\Phi^{-1}$)
``turned upside down''.
It is not hard to check that (\ref{eq:longsum})
can be replaced by the recursive definition
\begin{equation}
\label{eq:shortsum}
(\alpha_3 f)(x) = f(x) + 
\max \{ (\alpha_3 f)(y): \ y \in \widehat{P}, \ y \covers x \}
\end{equation}
which turns out to be the form most suitable for lifting
to the birational setting.

\begin{theorem}
\label{thm:pl-three}
$\rowP = \alpha_1 \circ \alpha_3 \circ \alpha_2$.
\end{theorem}

(Note that $\alpha_1 \circ \rowP \circ \alpha_1 =
\alpha_3 \circ \alpha_2 \circ \alpha_1$,
as in the original definition of $\row$.)

Similarly, in the birational setting put
\begin{eqnarray}
(\alpha_1 f)(x) & = & 1/f(x), \label{eq:alph1}\\
(\alpha_2 f)(x) & = & 
\bigsig^{\psum} \{ f(x)/f(y): \ y \in x^- \}, \ {\rm and} \label{eq:alph2}\\
(\alpha_3 f)(x) & = & 
f(x) \ \bigsig^{+} \{ (\alpha_3 f)(y): \ y \in x^+ \}\label{eq:alph3},
\end{eqnarray}
where $x^+$ denotes $\{ y \in \widehat{P}: \ y \covers x \}$
and $x^-$ denotes $\{ y \in \widehat{P}: \ x \covers y \}$.
(Note that definition (\ref{eq:alph3}),
like definition (\ref{eq:shortsum}), is recursive.)

\begin{theorem}
\label{thm:birational-three}
$\rowB = \alpha_1 \circ \alpha_3 \circ \alpha_2$.
\end{theorem}
Of course the $\alpha$'s in Theorem~\ref{thm:birational-three}
are not the $\alpha$'s in Theorem~\ref{thm:pl-three}
but their birational counterparts.

In the full paper, we derive Theorem~\ref{thm:pl-three}
from Theorem~\ref{thm:birational-three}
by tropicalization and dehomogenization.

\end{section}

\begin{section}{Recombination and Reciprocal Symmetry}
\label{sec:recombine}

As was noted by~\cite{prorow},
there is an intimate relationship between rowmotion and promotion
in rc-embedded posets: the two maps have the same orbit structure
because they are conjugate as elements of the toggle group.
This relationship becomes even clearer in 
the piecewise-linear and birational settings.
Let $P = [2] \times [2]$.
Here is the $\rowB$-orbit of $(1,2,3,4)$:
\begin{displaymath}
\begin{array}{ccccccccc}
( &  1  & , &  2   & , & 3    & , & 4   & )\\[5pt]
( & 1/4 & , & 5/8  & , & 5/12 & , & 5/4 & )\\[5pt]
( & 4/5 & , & 1/3  & , & 1/2  & , & 5/6 & )\\[5pt]
( & 6/5 & , & 12/5 & , & 8/5  & , &  1  & )
\end{array}
\end{displaymath}
Here is the $\proB$-orbit of $(1,2,5/12,5/4)$:
\begin{displaymath}
\begin{array}{ccccccccc}
( &  1  & , &  2   & , & 5/12 & , & 5/4 & )\\[5pt]
( & 1/4 & , & 5/8  & , & 1/2  & , & 5/6 & )\\[5pt]
( & 4/5 & , & 1/3  & , & 8/5  & , &  1  & )\\[5pt]
( & 6/5 & , & 12/5 & , &  3   & , &  4  & )
\end{array}
\end{displaymath}
Note that the same numbers appear as entries in both orbits,
with the same multiplicity.
More specifically, given $P = [a] \times [b]$,
define the {\em recombination map} $D$
as the map from the set of $P$-arrays to itself
such that for every $P$-array $f$,
the $(i,j)$ entry in $D(f)$ is the $(i,j)$ entry in $\rowB^{i-1}(f)$.

\begin{theorem}
\label{thm:recombine}
{\rm (the ``recombination lemma''):}
$D \circ \proB = \rowB \circ D$.
\end{theorem}

It follows from Theorem~\ref{thm:recombine} that $D$ is invertible
and that $\proB$ and $\rowB$ have the same orbit-structure.

A seemingly much deeper fact is the following consequence
of the work of~\cite{grinbergroby} (Theorem 10.6 in particular).

\begin{theorem}
\label{thm:recip}
{\rm (reciprocal symmetry):}
The $(a-i+1,b-j+1)$ entry in $\rowB^{a+b+1-i-j}(f)$ 
is the reciprocal of the $(i,j)$ entry in $f$.
\end{theorem}

Applying this theorem twice yields the conclusion
that for $n=(a+b+1-i-j)+(a+b+1-(a-i+1)-(b-j+1))=a+b$,
the $(i,j)$ entry in $\rowB^{n}$
is the reciprocal of the reciprocal of the $(i,j)$ entry in $f$.
This implies that $\rowB^{n}$ is the identity map
(and recombination then assures us that
that $\proB^{n}$ is the identity map as well).
The fact that $\row^{n}$ acts trivially on $J([a] \times [b])$
was first proved by~\cite{fon1993orbits}.

These facts have implications in the piecewise-linear setting.
The recombination property says that
the $(i,j)$ entry in $D(f)$ is the $(i,j)$ entry in $\rowB^{i-1}(f)$,
and reciprocal symmetry says that
the $(a-i+1,b-j+1)$ entry in $\rowP^{a+b+1-i-j}(f)$ 
is 1 minus the $(i,j)$ entry in $f$.
We also may conclude that $\rowP^{n}$ and $\proP^{n}$
are the identity map.
The last of these conclusions,
in combination with our remarks in section~\ref{sec:PL}
linking certian $P$-arrays with semistandard Young tableaux,
gives us a new proof of the standard fact that 
Sch{\"u}tzenberger promotion
on standard tableaux of fixed rectangular shape
with entries bounded by $n$ has order $n$.

We stress that recombination is not specific to $[a] \times [b]$,
but applies to any rc-embedded poset,
even in cases where rowmotion is not of finite order.
The recombination lemma is heavily based on Theorem 5.4 in~\cite{prorow}
(construction of an equivariant bijection).

\end{section}

\begin{section}{File-toggling and promotion}
\label{sec:filetoggle}

Here we restrict to $P$ of the form $[a] \times [b]$, with $n=a+b$.
The birational toggle operations $\tau_i$,
combined in unconstrained fashion,
generate a group that is infinite when $a>1$ or $b>1$
(we prove this in detail in the full article for the case $a=b=2$),
and its structure is likely to be quite complicated,
but some of the subgroups
admit homomorphisms to the symmetric group $S_n$,
and they can be useful for understanding rowmotion and promotion.
One such subgroup, generated by $n-1$ involutions
associated with the respective ranks of $P$,
was discovered by~\cite{grinbergroby}.
Here we study a different subgroup,
generated by $n-1$ involutions
associated with the respective files of $P$.

Recall that $[a] \times [b]$ can be partitioned into files
numbered 1 through $n-1$ from left to right.
Given $f: \widehat{P} \rightarrow \R^+$
with $f(\widehat{0}) = f(\widehat{1}) = 1$,
let $p_i$ ($1 \leq i \leq n-1$) be the product of the numbers $f(x)$
with $x$ belonging to the $i$th file of $P$,
let $p_0 = p_n = 1$,
and for $1 \leq i \leq n$ let $q_i = p_i/p_{i-1}$.
Call $q_1,\dots,q_n$ the {\em quotient sequence} associated with $f$,
and denote it by $Q(f)$.
This is analogous to the difference sequence introduced in~\cite{propproby}.
Note that the product $q_1 \cdots q_n$ telescopes to $p_n/p_0=1$.
For $i$ between 1 and $n-1$,
let $\tau_i^*$ be the product of the commuting involutions $\tau_x$
for all $x$ belonging to the $i$th file.
Lastly, given a sequence of $n$ numbers $w = (w_1,\dots,w_n)$,
and given $1 \leq i \leq n-1$, define 
$\sigma_i (w) = (w_1,\dots,w_{i-1},w_{i+1},w_{i},w_{i+2},\dots,w_n)$;
that is, $\sigma_i$ switches the $i$th and $i+1$st entries of $w$.

\begin{lemma}
\label{lem:swap}
For all $1 \leq i \leq n-1$, and for all $f$,
\begin{displaymath} Q(\tau_i^* f) = \sigma_i Q(f). \end{displaymath}
That is, toggling the $i$th file of $f$
swaps the $i$th and $i+1$st entries
of the quotient sequence of $f$.
\end{lemma}

Recalling that $\proB$ is the composition
$\tau_{n-1}^* \circ \cdots \circ \tau_{1}^*$,
we have:

\begin{corollary}
\label{cor:shift}
$Q(\proB f)$ is the leftward cyclic shift of $Q(f)$.
\end{corollary}

\end{section}

\begin{section}{Homomesy}
\label{sec:homomesy}

Given a set $X$, an operation $T : X \rightarrow X$
whose $n$th power is the identity map on $X$,
and a function $F$ from $X$ to a field $\cK$ of characteristic 0, 
we say that $F$ is {\em homomesic} 
relative to (or under the action of) $T$, 
or that the triple $(X,T,F)$ exhibits {\em homomesy}, 
if for all $x \in X$ the average
\begin{displaymath} \frac{1}{n} \sum_{k=0}^{n-1} F(T^k(x)) \end{displaymath}
equals some $c$ independent of $x$. 
We also say in this situation that the function $F$ 
(which we will sometimes call a {\em functional} on $X$)
is $c$-{\em mesic} relative to the map $T$.
The article by~\cite{propproby} gives examples of combinatorial situations 
in which homomesy holds.
See also~\cite{bloom2013homomesy}.

Theorem~\ref{thm:recip} yields as a corollary
that $((\R^+)^P,\rowB,F)$ is 0-mesic,
where $F(f) = \log(f(i,j) f(a+1-i,b+1-j))$
(factors cancel in pairs).
Applying recombination, we see that the same is true
if rowmotion is replaced by promotion.
In both cases, tropicalizing
yields homomesy for $F(f) = f(i,j) + f(a+1-i,b+1-j)$
under piecewise-linear rowmotion and promotion.

A different sort of homomesy comes from the files of $[a] \times [b]$.
Using Corollary~\ref{cor:shift},
one can show that for each $i$ between 1 and $n-1$,
if one defines $F_i(f)$ as the logarithm of 
the product of the values of $f(x)$
as $x$ ranges over the $i$th file of $[a] \times [b]$,
then $((\R^+)^P,\rowB,F_i)$ is 0-mesic.
This can be carried to the piecewise-linear setting as well.
Restricting to the vertices of $\cO(P)$,
one obtains the main homomesy theorem of~\cite{propproby}.

We can see both forms of homomesy on display in the rowmotion orbit
shown at the start of section~\ref{sec:recombine}.
For instance, the middle file of the poset
consists of the elements $w$ and $z$,
associated with the entries $v_1$ and $v_4$
of each vector $v$.
Defining $F(f)$ as $\log f(w)f(z)$,
we see that over the orbit the function $F$
takes on the values $\log 4$, $\log 5/16$, $\log 2/3$, and $\log 6/5$,
which sum to 0.

\begin{theorem}
\label{pl-homomesy}
Given $P = [a] \times [b]$, with $n=a+b$,
define functionals 
$F_{i,j}$ ($1 \leq i \leq a$, $1 \leq j \leq b$)
and $F_k$ ($1 \leq k \leq n-1$) by 
\begin{displaymath} F_{i,j}(f) = f(i,j) + f(a+1-i,b+1-j)), \end{displaymath}
\begin{displaymath} F_k(f) = \sum_{j-i=k-a} f(i,j). \end{displaymath}
These functionals are all homomesic
under the action of $\rowP$ and $\proP$.
\end{theorem}

\begin{theorem}
\label{birational-homomesy}
Given $P = [a] \times [b]$, with $n=a+b$,
define functionals 
$F_{i,j}$ ($1 \leq i \leq a$, $1 \leq j \leq b$)
and $F_k$ ($1 \leq k \leq n-1$) by 
\begin{displaymath} F_{i,j}(f) = \log(f(i,j) f(a+1-i,b+1-j)), \end{displaymath}
\begin{displaymath} F_k(f) = \log(\prod_{j-i=k-a} f(i,j)). \end{displaymath}
These functionals are all homomesic
under the action of $\rowB$ and $\proB$.
\end{theorem}

The recombination lemma easily implies
that a functional $F$ is homomesic under rowmotion
if and only if it is homomesic under promotion.
Also, any linear combination of homomesic functions is homomesic.

In the full version of the article,
a kind of converse of Theorem~\ref{pl-homomesy} will be proved:

\begin{theorem}
\label{pl-converse}
Given $P = [a] \times [b]$, with $p=ab$,
let $F$ be some function in the span of
the $p$ evaluation functions $f \mapsto f(i,j)$
(with $1 \leq i \leq a$, $1 \leq j \leq b$),
such that $F$ is homomesic under the action
of $\rowP$ (or equivalently under the action of $\proP$);
then $F$ must be a linear combination
of the functional $F_{i,j}$ and $F_k$ defined
in Theorem~\ref{pl-homomesy}.
\end{theorem}

Let $V$ be the vector space
spanned by the functionals $F_{i,j}$ and $F_k$.
It should be noted that the functionals $F_{i,j}$ and $F_k$
have linear dependencies,
so although they span $V$,
they are not a basis of $V$.

Although we have restricted ourselves to $(\R^+)^P$ for simplicity,
to the extent that our main results are
complicated but finite subtraction-free identities,
results like these homomesy theorems,
or the fact that rowmotion and promotion are of order $n$,
apply throughout the complement of some proper subvariety of $\C^P$
(though we need to use $\log |z|$ in place of $\log z$).
Also note that our birational maps are homogeneous,
so projective counterparts of rowmotion and promotion can be defined
and are likely to be helpful.

\end{section}

\acknowledgements
\label{sec:ack}
This work was supported by a grant from NSF.
The authors are grateful to Arkady Berenstein, Darij Grinberg, 
Alex Postnikov, Tom Roby, Richard Stanley, and Jessica Striker 
for helpful conversations, and to the referees for helpful suggestions.

\nocite{*}
\bibliographystyle{abbrvnat}
\bibliography{fpsac}
\label{sec:biblio}

\end{document}